\documentclass[review]{siamart190516}   

\usepackage{siampaper}
\newcommand{\USYMQR}{\textsc{Usymqr}\xspace}
\newcommand{\USYMLQ}{\textsc{Usymlq}\xspace}
\newcommand{\USYMLQR}{\textsc{Usymlqr}\xspace}

\newcommand{\BiLQ}{\textsc{BiLQ}\xspace}
\newcommand{\BiCG}{\textsc{BiCG}\xspace}
\newcommand{\BiCGSTAB}{\textsc{BiCGStab}\xspace}
\newcommand{\CGS}{\textsc{Cgs}\xspace}
\newcommand{\QMR}{\textsc{Qmr}\xspace}

\newcommand{\BiLQR}{\textsc{BiLQR}\xspace}
\newcommand{\TriLQR}{\textsc{TriLQR}\xspace}

\newcommand{\CGM}{\textsc{Cg}\xspace}
\newcommand{\MINRES}{\textsc{Minres}\xspace}
\newcommand{\SYMMLQ}{\textsc{Symmlq}\xspace}
\newcommand{\MINRESQLP}{\textsc{Minres-qlp}\xspace}

\newcommand{\LSQR}{\textsc{Lsqr}\xspace}

\newcommand{\LSLQ}{\textsc{Lslq}\xspace}
\newcommand{\LNLQ}{\textsc{Lnlq}\xspace}

\newcommand{\supL}{^{\textup{\tiny L}}}
\newcommand{\supC}{^{\textup{\tiny C}}}
\newcommand{\supQ}{^{\textup{\tiny Q}}}
\newcommand{\supLQ}{^{\textup{\tiny LQ}}}
\newcommand{\supCG}{^{\textup{\tiny CG}}}
\newcommand{\supQR}{^{\textup{\tiny QR}}}

\newcommand{\highlightline}[1]{%
  \hspace*{-\fboxsep}%
  \colorbox{gray!10}{%
    \parbox{\linewidth - 2\fboxsep - 0.25ex}{%
      #1
    }
  }
}
\newlength{\forwidth}
\newcommand{\highlightlinefor}[1]{%
  \hspace*{-\fboxsep}%
  \colorbox{gray!10}{%
    \parbox{\linewidth - 2\fboxsep - 0.5ex - \forwidth}{%
      #1
    }
  }
}
\newcommand{\highlight}[1]{%
  \hspace*{-\fboxsep}\colorbox{gray!10}{#1}%
}

\usepackage{tikzscale}
\usetikzlibrary{external}
\usepgfplotslibrary{colormaps}

\newcommand*{\includetikzgraphics}[2][]{%
  \includegraphics[#1]{#2}
}

\makeatletter
\renewcommand{\todo}[2][]{\tikzexternaldisable\@todo[#1]{#2}\tikzexternalenable}
\makeatother

\usepackage{newfloat}
\DeclareFloatingEnvironment[
  name=DoubleAlgorithm,
  placement=ht,
  fileext=lod]{doublealgorithm}

\newcommand{\diag}{\mathop{\text{diag}}}
\newcommand{\vol}{\mathop{\text{vol}}}

\crefname{subsection}{section}{sections}

\newcommand{\cahiernumber}{71}  

\usepackage{kbordermatrix}

\title{%
  {\BiLQ}:\@ An Iterative Method for Nonsymmetric Linear Systems with a Quasi-Minimum Error Property
}

\author{%
  Alexis Montoison%
  \thanks{%
    GERAD and Department of Mathematics and Industrial Engineering,
    Polytechnique Montr\'eal, QC, Canada.
    E-mail: \mailto{alexis.montoison@polymtl.ca}.
    Research partially supported by a merit scholarship of the Arbour foundation.
  }
  \and
  Dominique Orban%
  \thanks{%
    GERAD and Department of Mathematics and Industrial Engineering,
    Polytechnique Montr\'eal, QC, Canada.
    E-mail: \mailto{dominique.orban@gerad.ca}.
    Research partially supported by an NSERC Discovery Grant.
  }
}
\date{\today}

\begin{document}

  \nolinenumbers%
  \maketitle

  \thispagestyle{firstpage}
  \pagestyle{myheadings}

  \begin{abstract}
    We introduce an iterative method named \BiLQ for solving general square linear systems $Ax=b$ based on the Lanczos biorthogonalization process defined by least-norm subproblems, and that is a natural companion to \BiCG and \QMR.
    Whereas the \BiCG \citep{fletcher-1976}, \CGS \citep{sonneveld-1989} and \BiCGSTAB \citep{van_der_vost-1992} iterates may not exist when the tridiagonal projection of $A$ is singular, \BiLQ is reliable on compatible systems even if A is ill-conditioned or rank deficient.
    As in the symmetric case, the \BiCG residual is often smaller than the \BiLQ residual and, when the \BiCG iterate exists, an inexpensive transfer from the \BiLQ iterate is possible.
    Although the Euclidean norm of the \BiLQ error is usually not monotonic, it is monotonic in a different norm that depends on the Lanczos vectors.
    We establish a  similar property for the \QMR \citep{freund-nachtigal-1991} residual.
    \BiLQ combines with \QMR to take advantage of two initial vectors and solve a system and an adjoint system simultaneously at a cost similar to that of applying either method.
    We derive an analogous combination of \USYMLQ and \USYMQR based on the orthogonal tridiagonalization process \citep{saunders-simon-yip-1988}.
    The resulting combinations, named \BiLQR and \TriLQR, may be used to estimate integral functionals involving the solution of a primal and an adjoint system.
    We compare \BiLQR and \TriLQR with \MINRESQLP on a related augmented system, which performs a comparable amount of work and requires comparable storage.
    In our experiments, \BiLQR terminates earlier than \TriLQR and \MINRESQLP in terms of residual and error of the primal and adjoint systems.
  \end{abstract}



  \section{Introduction}

  We consider the square consistent linear system
  \begin{equation}
    \label{eq:Ax=b}
    Ax = b,
  \end{equation}
  where \(A \in \R^{n \times n}\) can be nonsymmetric, is either large and sparse, or is only available as a linear operator, i.e., via operator-vector products.
  We assume that \(A\) is nonsingular.
  Systems such as~\eqref{eq:Ax=b} arise in the discretization of partial differential equations (PDEs) in numerous applications, including compressible turbulent fluid flow \citep{chisholm-zingg-2009}, and in circuit simulation \citep{davis-natarajan-2012}.
  We consider Krylov subspace methods and are interested in generating iterates with guarantees as to the decrease of the error \(x_k - x_\star\) in a certain norm, where \(x_\star\) is the solution of~\eqref{eq:Ax=b}.

  The foundation of Krylov methods is a basis-generation process upon which three methods may be developped: one computing the minimum-norm solution of an under-determined system, one solving a square system and imposing a Galerkin condition, and one solving an over-determined system in the least-squares sense.
  These methods may be implemented with the help of a LQ, LU or QR factorization of a related operator, respectively.

  In this paper, we develop an iterative method named \BiLQ of the first type based on the~\cite{lanczos-1950} biorthogonalization process.
  Together with \BiCG \citep{fletcher-1976} and \QMR \citep{freund-nachtigal-1991}, \BiLQ completes the family of methods based on the biorthogonalization process.
  We begin by stating the defining properties of \BiLQ, describing its implementation in detail, and illustrating its behavior on numerical examples side by side with \BiCG and \QMR.

  In a second stage, we exploit the fact that the biorthogonalization process requires two initial vectors to develop a combination of \BiLQ and \QMR that solves~\eqref{eq:Ax=b} together with a dual system
  \begin{equation}
    \label{eq:A'y=c}
    A^T t = c
  \end{equation}
  simultaneously at a cost comparable to that of applying \BiLQ or \QMR only to solve one of those systems.
  The resulting combination is named \BiLQR and is employed to illustrate the computation of superconvergent estimates of integral functionals arising in certain PDE problems.

  We note that a similar approach may be developed for the \cite{saunders-simon-yip-1988} orthogonal tridiagonalization process, which also requires two initial vectors, by combining \USYMLQ and \USYMQR.
  The resulting combination is named \TriLQR.

  Finally, we compare \BiLQR and \TriLQR with \MINRESQLP on a related augmented system to solve both~\eqref{eq:Ax=b} and~\eqref{eq:A'y=c} simultaneously.
    In our experiments, \BiLQR terminates earlier than \TriLQR and \MINRESQLP in terms of residual and error of the primal and adjoint systems.

  Our Julia \citep{julia-2017} implementation of \BiLQ, \QMR, \USYMLQ, \USYMQR, \BiLQR, \TriLQR, and \MINRESQLP are available from \https{github.com/JuliaSmoothOptimizers/Krylov.jl}.
  Thanks to multiple dispatch, a language feature allowing automatic compilation of variants of each method corresponding to inputs expressed in various floating-point systems, our implementations run in any floating-point precision supported.

  \subsection*{Related Research}

    \cite{paige-saunders-1975} develop one of the best-known minimum error methods, \SYMMLQ, based on the symmetric Lanczos process.
    \SYMMLQ inspires \cite{estrin-orban-saunders-2019a,estrin-orban-saunders-2019b} to develop \LSLQ and \LNLQ for rectangular problems based on the \cite{golub-kahan-1965} process.
    \LSLQ and \LNLQ are equivalent to \SYMMLQ applied to the normal equations and normal equations of the second kind, respectively.

    \cite{saunders-simon-yip-1988} define \USYMLQ for square consistent systems based on the orthogonal tridiagonalization process.
    \USYMLQ is based on a subproblem similar to that of \SYMMLQ, and coincides with \SYMMLQ in the symmetric case.
    Its companion method, \USYMQR, is similar in spirit to \MINRES.
    \cite{buttari-orban-ruiz-titley_peloquin-2019} combine both into a method named \USYMLQR designed to solve symmetric saddle-point systems with general right-hand side, and inspire the developement of \BiLQR and \TriLQR in the present paper.

    \cite{weiss-1994} decribes two types of error-minimizing Krylov methods for square \(A\); one based on a process applied to \(A^T A\), and one to \(A^T\).
    Our approach is to apply the biorthogonalization process directly to \(A\).
    We defer a numerical stability analysis to future work, but note that \cite{paige-panayotov-zemke-2014} study the augmented stability of the biorthogonalization process.
    In this sense, we make the implicit assumption that computations are carried out in exact arithmetic.
    This assumption prompted us to develop our implementations so that they can be applied in any supported floating-point arithmetic.

    The simultaneous solution of a system and an adjoint system has attracted attention in the past.
    Notably, \cite{lu-darmofal-2003} devise a variant of \QMR to solve both systems at once at a cost approximately equal to that of \QMR applied to one of the systems but with an increase in storage requirements.
    \cite{golub-stoll-wathen-2008} follow a similar approach and use a variant of \USYMQR to solve both~\eqref{eq:Ax=b} and~\eqref{eq:A'y=c}.
    An advantage of \USYMQR is to produce monotonic residuals in the Euclidean norm for both systems.
    We illustrate in \cref{tab:cost_bilqr_trilqr_minresqlp} that our methods are cheaper and have smaller storage requirements than those of \cite{lu-darmofal-2003} and \cite{golub-stoll-wathen-2008} though residuals are not monotonic in the Euclidean norm.

  \subsection*{Notation}

  Matrices and vectors are denoted by capital and lowercase Latin letters, respectively, and scalars by Greek letters.
  An exception is made for Givens cosines and sines $(c, s)$ that compose reflections.
  For a vector v, $\|v\|$ denotes the Euclidean norm of $v$, and for symmetric and positive-definite \(N\), the \(N\)-norm of \(v\) is \(\|v\|_N^2 = v^T N v\).
  For a matrix $M$, $\|M\|_F$ denotes the Frobenius norm of $M$.
  The vector $e_i$ is the $i$-th column of an identity matrix of size dictated by the context.
  Vectors and scalars decorated by a bar will be updated at the next iteration.
  For \(j = 2, \dots, k\), we use the compact representation
  \[
    Q_{j-1,j}
    =
    \kbordermatrix{
      & j-1 & j               \\
      & c_j & \phantom{-} s_j \\
      & s_j &          -  c_j
    }
    :=
    \begin{bmatrix}
      I_{j-2} &     &                & \\
              & c_j & \phantom{-}s_j & \\
              & s_j & -c_j           & \\
              &     &                & I_{k-j}
    \end{bmatrix},
  \]
  for orthogonal reflections, where \(s_j^2 + c_j^2 = 1\), where border indices indicate row and column numbers, and where \(I_k\) represents the \(k\)\(\times\)\(k\) identity operator.
  We abuse the notation $\bar{z}_k = (z_{k-1},~\bar{\zeta}_k)$ to represent the column vector $\begin{bmatrix} z_{k-1}^T & \bar{\zeta}_k \end{bmatrix}^T$.

  \section{Derivation of \BiLQ}

  \subsection{The Lanczos Biorthogonalization Process}

  The \citeauthor{lanczos-1950} biorthogonalization process generates sequences of vectors \(\{v_k\}\) and \(\{u_k\}\) such that \(v_i^T u_j = \delta_{ij}\) in exact arithmetic for as long as the process does not break down.
  The process is summarized as \autoref{alg:bilanczos}.

  \begin{algorithm}[H]
    \caption{%
      Lanczos Biorthogonalization Process
    }
    \label{alg:bilanczos}
    \begin{algorithmic}[1]
      \Require $A$, $b$, $c$
      \State $v_0 = 0$, $u_0 = 0$
      \State $\beta_1 v_1 = b$, $\gamma_1 u_1 = c$ \Comment{$(\beta_1, \gamma_1)$ so that $v_1^T u_1 = 1$}
      \For{$k$ = 1,~2,~\(\dots\)}
        \State $q = A v_k - \gamma_k v_{k-1}$, $\alpha_k = u_k^T q$
        \State $p = A^T u_k - \beta_k u_{k-1}$
        \State $\beta_{k+1} v_{k+1} = q  - \alpha_k v_k$ \Comment{$(\beta_{k+1}, \gamma_{k+1})$ so that $v_{k+1}^T u_{k+1} = 1$}
        \State $\gamma_{k+1} u_{k+1} = p - \alpha_k u_k$
      \EndFor
    \end{algorithmic}
  \end{algorithm}

  We denote \(V_k = \begin{bmatrix} v_1 & \dots & v_k \end{bmatrix}\) and \(U_k = \begin{bmatrix} u_1 & \dots & u_k \end{bmatrix}\).
  Without loss of generality, we choose the scaling factors \(\beta_k\) and \(\gamma_k\) so that \(v_k^T u_k = 1\) for all \(k \geq 1\), i.e., \(V_k^T U_k = I_k\).
  After \(k\) iterations, the situation may be summarized as
  \begin{subequations}
    \label{eq:bilanczos}
    \begin{align}
         A V_k & = V_k T_k~ + \beta_{k+1} v_{k+1} e_k^T = V_{k+1} T_{k+1,k}
         \label{eq:bilanczos-V}
      \\ A^T U_k & = U_k T_k^T + \gamma_{k+1} u_{k+1} e_k^T = U_{k+1} T_{k,k+1}^T,
        \label{eq:bilanczos-U}
    \end{align}
  \end{subequations}
  where
  \begin{equation*}
    T_k =
    \begin{bmatrix}
      \alpha_1 & \gamma_2 &          & \\
      \beta_2  & \alpha_2 & \ddots   & \\
               & \ddots   & \ddots   & \gamma_k \\
               &          & \beta_k  & \alpha_k
    \end{bmatrix},
    ~T_{k,k+1} =
    \begin{bmatrix}
      T_{k} & \gamma_{k+1} e_k
    \end{bmatrix},
    ~T_{k+1,k} =
    \begin{bmatrix}
      T_{k} \\
      \beta_{k+1} e_{k}^T
    \end{bmatrix}.
  \end{equation*}

  The columns of $V_k$ and $U_k$ form a basis for \(\mathcal{K}_k := \Span\{b, Ab, \cdots, A^{k-1}b\}\) and \(\mathcal{L}_k := \Span\{c, A^T c, \cdots, (A^T)^{k-1}c\}\), respectively.
  Though \(V_k\) cannot be expected to be orthogonal to \(U_k\) in inexact arithmetic, and therefore \(U_k^T A V_k = T_k\) cannot be expected to hold, \eqref{eq:bilanczos} usually holds to within machine precision.

  \subsection{Definition of \BiLQ}

  By definition, \BiLQ generates an approximation \(x\supL_k\) to a solution of~\eqref{eq:Ax=b} of the form \(x\supL_k = V_k y\supL_k\), where \(y\supL_k \in \R^k\) solves
  \begin{equation}
    \label{eq:bilq-sub}
    \minimize{y} \ \|y\| \quad \st \ T_{k-1,k} y = \beta_1 e_1.
  \end{equation}
  By contrast, \BiCG \citep{fletcher-1976} generates \(x\supC_k = V_k y\supC_k\) where \(y\supC_k \in \R^k\) solves
  \begin{equation}
    \label{eq:bicg-sub}
    T_k y = \beta_1 e_1,
  \end{equation}
  and \QMR \citep{freund-nachtigal-1991} generates \(x\supQ_k = V_k y\supQ_k\) where \(y\supQ_k \in \R^k\) solves
  \begin{equation}
    \label{eq:qmr-sub}
    \minimize{y} \ \|T_{k+1,k} y - \beta_1 e_1\|.
  \end{equation}

  When A is symmetric and $b = c$, \Cref{alg:bilanczos} coincides with the symmetric Lanczos process and the three above methods are equivalent to \SYMMLQ \citep{paige-saunders-1975}, \CGM \citep{hestenes-stiefel-1952}, and \MINRES \citep{paige-saunders-1975}, respectively.

  \subsection{An LQ factorization}

  We determine \(y\supL_k\) solution to~\eqref{eq:bilq-sub} via the LQ factorization of \(T_{k-1,k}\), which we obtain from the LQ factorization
  \begin{subequations}
    \label{eq:lq_bicg}
    \begin{align}
      T_k & = \bar{L}_k Q_k, \text{ where}
      \\
      \overline{L}_k & =
      \begin{bmatrix}
        \delta_1      &           &                   &                   &               & \\
        \lambda_1     & \delta_2  &                   &                   &               & \\
        \varepsilon_1 & \lambda_2 & \delta_3          &                   &               & \\
                      & \ddots    & \ddots            & \ddots            &               & \\
                      &           & \varepsilon_{k-3} & \lambda_{k-2}     & \delta_{k-1}  & \\
                      &           &                   & \varepsilon_{k-2} & \lambda_{k-1} & \bar{\delta}_k
      \end{bmatrix}
      =
      \begin{bmatrix}
          L_{k-1}                                               & 0 \\
          \varepsilon_{k-2} e^T_{k-2} + \lambda_{k-1} e^T_{k-1} & \bar{\delta}_k
      \end{bmatrix},
    \end{align}
  \end{subequations}
  and $Q_k^T = Q_{1,2} Q_{2,3} \cdots Q_{k-1,k}$ is orthogonal and defined as a product of Givens reflections.
  Indeed, the above yields the LQ factorization
  \begin{equation}
    \label{eq:lq_bilq}
    T_{k-1,k}
    =
    \begin{bmatrix}
      L_{k-1} & 0
    \end{bmatrix}
    Q_k.
  \end{equation}

  If we initialize \(\bar{\delta}_1 := \alpha_1\), \(\bar{\lambda}_1 := \beta_2\), \(c_1 = -1\), and \(s_1 = 0\), individual factorization steps may be represented as an application of $Q_{k-2,k-1}$ to $T_k Q_{k-2}^T$:
  \[
    \kbordermatrix{
        & k-2                  & k-1          & k        \\
    k-2 & \bar{\delta}_{k-2}   & \gamma_{k-1} &          \\
    k-1 & \bar{\lambda}_{k-2}  & \alpha_{k-1} & \gamma_k \\
    k   &                      & \beta_k      & \alpha_k
    }
    \kbordermatrix{
    & k-2      & k-1                 & k \\
    & c_{k-1}  & \phantom{-} s_{k-1} &   \\
    & s_{k-1}  & -c_{k-1}            &   \\
    &          &                     & 1
    }
    =
    \kbordermatrix{
    & k-2               & k-1                 & k        \\
    & \delta_{k-2}      & 0                   &          \\
    & \lambda_{k-2}     & \bar{\delta}_{k-1}  & \gamma_k \\
    & \varepsilon_{k-2} & \bar{\lambda}_{k-1} & \alpha_k
    },
  \]
  followed by an application of $Q_{k-1,k}$ to the result:
  \[
    \kbordermatrix{
        & k-2               & k-1                 & k        \\
    k-2 & \delta_{k-2}      &                     &          \\
    k-1 & \lambda_{k-2}     & \bar{\delta}_{k-1}  & \gamma_k \\
    k   & \varepsilon_{k-2} & \bar{\lambda}_{k-1} & \alpha_k
    }
    \kbordermatrix{
    & k-2 & k-1 & k               \\
    & 1   &     &                 \\
    &     & c_k & \phantom{-} s_k \\
    &     & s_k & -c_k
    }
    =
    \kbordermatrix{
    & k-2               & k-1           & k              \\
    & \delta_{k-2}      &               &                \\
    & \lambda_{k-2}     & \delta_{k-1}  &                \\
    & \varepsilon_{k-2} & \lambda_{k-1} & \bar{\delta}_k
    }.
  \]
  The reflection $Q_{k-1,k}$ is designed to zero out $\gamma_k$ on the superdiagonal of $T_k$ and affects three rows and two colums.
  It is defined by
  \begin{equation}
    \delta_{k-1} = \sqrt{\bar{\delta}_{k-1}^2 + \gamma_k^2},
    \quad
    c_k = \bar{\delta}_{k-1} / \delta_{k-1},
    \quad
    s_k = \gamma_k / \delta_{k-1},
  \end{equation}
  and yields the recursion
    \begin{subequations}
      \label{lq_formula}
      \begin{alignat}{2}
        \varepsilon_{k-2}   &= s_{k-1} \beta_k, \quad & k & \ge 3, \\
        \bar{\lambda}_{k-1} &= -c_{k-1} \beta_k, \quad & k & \ge 3, \\
        \lambda_{k-1}       &= c_k \bar{\lambda}_{k-1} + s_k \alpha_k, \quad & k & \ge 2, \\
        \bar{\delta}_k      &= s_k \bar{\lambda}_{k-1} - c_k \alpha_k, \quad & k & \ge 2.
      \end{alignat}
    \end{subequations}

  \subsection{Definition and update of the \BiLQ and \BiCG iterates}

  In order to compute $y_k\supL$ solution of \eqref{eq:bilq-sub} using \eqref{eq:lq_bilq}, we solve \(\begin{bmatrix} L_{k-1} & 0 \end{bmatrix} Q_k y_k\supL = \beta_1 e_1\).
  If \(z_{k-1} := (\zeta_1, \dots, \zeta_{k-1})\) is defined so that $L_{k-1} z_{k-1} = \beta_1 e_1$, then the minimum-norm solution of~\eqref{eq:bilq-sub} is $y_k\supL = Q_k^T \begin{bmatrix} z_{k-1} \\ 0 \end{bmatrix}$, and $\|y_k\supL\| = \|z_{k-1}\|$.

  We may compute $y_k\supC$ in~\eqref{eq:bicg-sub} simultaneously as a cheap update of \(y_k\supL\).
  Indeed,~\eqref{eq:bicg-sub} and \eqref{eq:lq_bicg} yield $\overline{L}_k Q_k y_k\supC = \beta_1 e_1$.
  Let $\bar{z}_k := (z_{k-1}, \bar{\zeta}_k)$ be defined so $\overline{L}_k \bar{z}_k = \beta_1 e_1$.
  Then, $y_k\supC = Q_k^T \bar{z}_k$. If \(\bar{\delta}_k = 0\), \(y_k\supC\) and the \BiCG iterate \(x_k\supC\) are undefined.
  The components of \(\bar{z}_k\) are computed from
  \begin{subequations}
    \label{zk}
    \begin{align}
      \eta_k &=
      \begin{cases}
        \beta_1, & k = 1, \\
        -\lambda_1 \zeta_1, & k = 2, \\
        -\varepsilon_{k-2} \zeta_{k-2} - \lambda_{k-1} \zeta_{k-1}, & k \ge 3,
      \end{cases}
      \\
      \zeta_{k-1} & = \eta_{k-1} / \delta_{k-1}, \quad k \ge 2, \\
      \bar{\zeta}_k & = \eta_k / \bar{\delta}_k, \quad \text{if } \bar{\delta}_k \neq 0.
    \end{align}
  \end{subequations}

  By definition, $x_k\supL = V_k y_k\supL$ and $x_k\supC = V_k y_k\supC$. To avoid storing $V_k$, we let
  \begin{equation}
    \overline{D}_k := V_k Q_k^T = \begin{bmatrix} d_1,~  d_2,~ \cdots,~ d_{k-1},~ \bar{d}_k \end{bmatrix},~~~\bar{d}_1 = v_1,
  \end{equation}
  defined by the recursion
  \begin{equation}
    \begin{aligned}
      d_{k-1}   & = c_k \bar{d}_{k-1} + s_k v_k \\
      \bar{d}_k & = s_k \bar{d}_{k-1} - c_k v_k.
    \end{aligned}
  \end{equation}
  Finally,
  \begin{subequations}
    \label{eq:bilq_bicg}
    \begin{align}
      x_k\supL & = V_k y_k\supL = \overline{D}_k \begin{bmatrix} z_{k-1} \\ 0 \end{bmatrix} = D_{k-1} z_{k-1} = x_{k-1}\supL + \zeta_{k-1} d_{k-1}
      \label{eq:xL-update}
      \\
      x_k\supC & = V_k y_k\supC = \overline{D}_k \bar{z}_k = D_{k-1} z_{k-1} + \bar{\zeta}_k \bar{d}_k = x_k\supL + \bar{\zeta}_k \bar{d}_k.
      \label{eq:xC-transfer}
    \end{align}
  \end{subequations}

  We see from~\eqref{eq:xC-transfer} that it is possible to transfer from \(x_k\supL\) to \(x_k\supC\) cheaply provided \(\bar{\zeta}_k \neq 0\).
 Such transfer was described by \cite{paige-saunders-1975} as an inexpensive update from the \SYMMLQ to the \CGM point in the symmetric case.

\subsection{Residuals estimates}

The identity~\eqref{eq:bilanczos-V} allows us to write the residual associated to \(x_k = V_k y_k\) as
\[
  r_k =
  b - A x_k = 
  \beta_1 v_1 - A V_k y_k =
  \beta_1 v_1 - V_{k+1} T_{k+1,k} y_k.
\]
Thus,~\eqref{eq:bilq-sub} yields the residual at the \BiLQ iterate:
  \begin{align}
    r_k\supL & = V_{k-1}(\beta_1 e_1 - T_{k-1,k} y_k\supL) - (\beta_k e_{k-1} + \alpha_k e_k)^T y_k\supL \, v_k - \beta_{k+1} e_k^T y_k\supL \, v_{k+1}
    \nonumber \\
          & = - (\beta_k e_{k-1} + \alpha_k e_k)^T y_k\supL \, v_k - \beta_{k+1} e_k^T y_k\supL \, v_{k+1},
    \label{eq:rkL}
  \end{align}
  and~\eqref{eq:bicg-sub} yields the residual at the \BiCG iterate:
  \[
    r_k\supC = V_k(\beta_1 e_1 - T_k y_k\supC) - \beta_{k+1} v_{k+1} e_k^T y_k\supC
          = - \beta_{k+1} e_k^T y_k\supC v_{k+1}.
  \]
  Because $Q_k^T = Q_{1,2} Q_{2,3} \cdots Q_{k-1,k}$, we have
  \begin{equation*}
    \begin{aligned}
      e_{k-1}^T Q_k^T & =
      e_{k-1}^T Q_{k-2,k-1} Q_{k-1,k} =
      s_{k-1} e_{k-2}^T - c_{k-1} c_{k} e_{k-1}^T - c_{k-1} s_k e_k^T,
      \\
      e_{k}^T Q_k^T & =
      e_k^T Q_{k-1,k} =
      s_{k} e_{k-1}^T - c_{k} e_k^T,
    \end{aligned}
  \end{equation*}
  so that
  \begin{subequations}
    \label{eq:residuals}
    \begin{align*}
      e_{k-1}^T y_k\supL & = e_{k-1}^T Q_k^T \begin{bmatrix} z_{k-1} \\ 0 \end{bmatrix} = s_{k-1} \zeta_{k-2} - c_{k-1} c_k \zeta_{k-1}, \\
      e_k^T y_k\supL & = e_k^T Q_k^T \begin{bmatrix} z_{k-1} \\ 0 \end{bmatrix} = s_k \zeta_{k-1}, \\
      e_k^T y_k\supC & = e_k^T Q_k^T \bar{z}_k = s_k \zeta_{k-1} - c_k \bar{\zeta}_k.
    \end{align*}
  \end{subequations}
  Therefore, if we define
  \(
    \mu_k =
    \beta_k (s_{k-1} \zeta_{k-2} - c_{k-1} c_k \zeta_{k-1}) + \alpha_k s_k \zeta_{k-1}
  \),
  \(
    \omega_k = \beta_{k+1} s_k \zeta_{k-1}
  \)
  and
  \(
    \rho_k = \beta_{k+1} (s_k \zeta_{k-1} - c_k \bar{\zeta}_k)
  \),
  we obtain
  \[
    \| r_k\supL \| = \sqrt{\mu_k^2 \|v_k\|^2 + \omega_k^2 \|v_{k+1}\|^2 + 2 \mu_k \omega_k v_k^T v_{k+1}},
  \]
  and
  \[
    \|r_k\supC \| = |\rho_k| \, \|v_{k+1}\|.
  \]

  We summarize the complete procedure as \autoref{alg:bilq}.
  For simplicity, we do not include a lookahead procedure, although a robust implementation should in order to avoid serious breakdowns \citep{parlett-taylor-liu-1985}.
  Table~\ref{tab:cost_bilq_bicg_qmr} summarizes the cost per iteration of \BiLQ, \BiCG and \QMR.
  Each method requires one operator-vector product with \(A\) and one with \(A^T\) per iteration.
  We assume that in-place ``gemv'' updates of the form $y \leftarrow Av + \gamma y$ and $y \leftarrow A^T u + \beta y$ are available.
  Otherwise, each method requires two additional $n$-vectors to store \(Av\) and \(A^T u\).
  In the table, ``dots'' refers to dot products of $n$-vectors, ``scal'' refers to scaling an $n$-vector by a scalar, and ``axpy'' refers to adding a multiple of one $n$-vector to another one.

  \begin{algorithm}[ht]
    \caption{%
      \BiLQ
    }
    \label{alg:bilq}
    \begin{algorithmic}[1]
      \Require $A$, $b$, $c$
      \State $\beta_1 v_1 = b$, $\gamma_1 u_1 = c$ \Comment{$(\beta_1, \gamma_1)$ so that $v_1^T u_1 = 1$}
      \State $\alpha_1 = u_1^T A v_1$ \Comment{begin biorthogonalization}
      \State $\beta_{2} v_{2} = A   v_1 - \alpha_1 v_1$
      \State $\gamma_{2} u_{2} = A^T u_1 - \alpha_1 u_1$
      \State $c_1 = -1$, $s_1 = 0$, $\bar{\delta}_1 = \alpha_1$ \Comment{begin $LQ$ factorization}
      \State $\eta_1 = \beta_1$, $\bar{d}_1 = v_1$, $x\supL_1 = 0$
      \For{$k$ = 2,~3,~\(\dots\)}
        \State $q = A v_k - \gamma_k v_{k-1}$, $\alpha_k = u_k^T q$ \Comment{continue biorthogonalization}
        \State $p = A^T u_k - \beta_k u_{k-1}$
        \State $\beta_{k+1} v_{k+1} = q  - \alpha_k v_k$ \Comment{$(\beta_{k+1}, \gamma_{k+1})$ so that $v_{k+1}^T u_{k+1} = 1$}
        \State $\gamma_{k+1} u_{k+1} = p - \alpha_k u_k$
        \State $\delta_{k-1} = (\bar{\delta}_{k-1}^2 + \gamma_k^2)^{\frac12}$ \Comment{compute $Q_{k-1,k}$}
        \State $c_k = \bar{\delta}_{k-1} / \delta_{k-1}$
        \State $s_k = \gamma_k / \delta_{k-1}$
        \State $\varepsilon_{k-2} = s_{k-1} \beta_k$ \Comment{continue $LQ$ factorization}
        \State $\lambda_{k-1} = - c_{k-1} c_k \beta_k + s_k \alpha_k$
        \State $\bar{\delta}_k = - c_{k-1} s_k \beta_k - c_k \alpha_k$
        \State $\zeta_{k-1} = \eta_{k-1} / \delta_{k-1}$ \Comment{update $z_{k-1}$}
        \State $\eta_k = -\varepsilon_{k-2} \zeta_{k-2} - \lambda_{k-1} \zeta_{k-1}$
        \State $\mu_k =\beta_k (s_{k-1} \zeta_{k-2} - c_{k-1} c_k \zeta_{k-1}) + \alpha_k s_k \zeta_{k-1}$
        \State $\omega_k = \beta_{k+1} s_k \zeta_{k-1}$
        \State $\| r\supL_k \| = (\mu_k^2 \|v_k\|^2 + \omega_k^2 \|v_{k+1}\|^2 + 2 \mu_k \omega_k v_k^T v_{k+1})^{\frac12}$ \Comment{compute $\| r\supL_k \|$}
        \If{$\bar{\delta}_k \neq 0$}
          \State $\bar{\zeta}_k = \eta_k / \bar{\delta}_k$ \Comment{optional: update $\bar{z}_k$}
          \State $\rho_k = \beta_{k+1} (s_k \zeta_{k-1} - c_k \bar{\zeta}_k)$
          \State $\|r\supC_k \| = |\rho_k| \, \|v_{k+1}\|$ \Comment{optional: compute $\| r\supC_k \|$}
        \EndIf
        \State $d_{k-1} = c_k \bar{d}_{k-1} + s_k v_k$ \Comment{update $\overline{D}_k$}
        \State $\bar{d}_k = s_k \bar{d}_{k-1} - c_k v_k$
        \State $x\supL_k = x\supL_{k-1} + \zeta_{k-1} d_{k-1}$ \Comment{\BiLQ point}
      \EndFor
      \If{$\bar{\delta}_k \neq 0$}
        \State $x\supC_k = x\supL_k + \bar{\zeta}_k \bar{d_k}$ \Comment{optional: \BiCG point}
      \EndIf
    \end{algorithmic}
  \end{algorithm}

  \begin{table}[H]
    \label{tab:cost_bilq_bicg_qmr}
    \caption{Storage and cost per iteration of methods based on \cref{alg:bilanczos}.}
    \begin{center}
      \setlength{\arrayrulewidth}{1pt}
      \begin{tabular}{lrrrr}
        \hline
        & $n$-vectors & dots & scal & axpy \\
        \hline
        \BiLQ & 6 & 2 & 3 & 7 \\
        \BiCG & 6 & 2 & 3 & 6 \\
        \QMR  & 7 & 2 & 4 & 7 \\
        \hline
      \end{tabular}
    \end{center}
  \end{table}

  \subsection{Properties}
  \label{sec:bilq-properties}

  By construction, assuming \cref{alg:bilanczos} does not break down, there exists an iteration \(p \leq n\) such that \(x\supL_{p+1} = x\supC_p = x_\star\), the exact solution of~\eqref{eq:Ax=b}.
  In particular, there exists \(y_\star\) such that \(x_\star = V_p y_\star\).

  The definition~\eqref{eq:bilq-sub} of \(y\supL_k\) ensures that \(\|y\supL_k\|\) is monotonically increasing while \(\|y\supL_k - y_\star\|\) is monotonically decreasing.
  Because \(V_k^T U_k = I_k\) at each iteration, the iteration-dependent norm
  \begin{equation}
    \label{eq:xLnorm}
    \|x\supL_k\|_{U_k U_k^T} = \|y\supL_k\|
  \end{equation}
  is monotonically increasing.
  Because we may write
  \begin{equation}
    \label{eq:xLk-Vp}
    x\supL_k =
    V_k y\supL_k =
    V_p
    \begin{bmatrix}
      y\supL_k \\ 0
    \end{bmatrix},
  \end{equation}
  \(\|x\supL_k\|_{U_p U_p^T} = \|x\supL_k\|_{U_k U_k^T}\) is also monotonically increasing, and the error norm
  \begin{equation}
    \label{eq:xLerror}
    \|x\supL_k - x_\star\|_{U_p U_p^T}
  \end{equation}
  is monotonically decreasing.
  Note that~\eqref{eq:xLnorm} is readily computable as \(\|z_{k-1}\|\), and can be updated as
  \[
    \|x\supL_{k+1}\|_{U_{k+1} U_{k+1}^T}^2 = \|x\supL_k\|^2_{U_k U_k^T} + \zeta_k^2.
  \]
  A lower bound on the error~\eqref{eq:xLerror} can be obtained as \(\|z_{k-d} - z_{k-1}\|\) for a user-defined \textit{delay} of \(d\) iterations.
  Such a lower bound may be used to define a simple, though not robust, error-based stopping criterion \citep{estrin-orban-saunders-2019b}.

  The following result establishes properties of \(x\supL_k\) that are analogous to those of the \SYMMLQ iterate in the symmetric case.

  \begin{shadyproposition}
    \label{prop:xL-min-norm}
    Let \(x_\star\) be as above.
    The \(k\)th \BiLQ iterate \(x\supL_k\) solves
    \begin{equation}
      \label{eq:xL-min-norm}
      \minimize{x} \, \|x\|_{U_k U_k^T} \ \st \, x \in \Range(V_k), \ b - A x \perp \Range(U_{k-1}),
    \end{equation}
    and
  \begin{equation}
    \label{eq:xL-min-error}
    \minimize{x} \, \|x - x_\star\|_{U_p U_p^T} \ \st \, x \in \Range(V_p V_p^T A^T U_{k-1}).
  \end{equation}
  \end{shadyproposition}

  \begin{proof}
    The first set of constraints of~\eqref{eq:xL-min-norm} imposes that there exist \(y \in \R^k\) such that \(x = V_k y\).
    By biorthogonality, the objective value at such an \(x\) can be written \(\|V_k y\|_{U_k U_k^T} = \|y\|\).
    Biorthogonality again and~\eqref{eq:rkL} show that \(y_k\) defined in~\eqref{eq:bilq-sub} is primal feasible for~\eqref{eq:xL-min-norm}.
    Dual feasibility of~\eqref{eq:xL-min-norm} requires that there exist a vector \(q\) such that \(y = V_k^T A^T U_{k-1} q\).
    By~\eqref{eq:bilanczos-U} and biorthogonality one more time, this amounts to \(y = T_{k-1,k}^T q\), which is the same as dual feasibility for~\eqref{eq:bilq-sub}.
    Thus, \(V_k y\supL_k\) is,  optimal for~\eqref{eq:xL-min-norm}.

    To establish primal feasibility of \(x\supL_k\) for~\eqref{eq:xL-min-error}, note first that~\eqref{eq:bilanczos-U} yields \(A^T U_{k-1} = U_k T_{k-1,k}^T\). 
    Let \(\bar{V}_{p-k}\) denote the last \(p-k\) columns of \(V_p\).
    Biorthogonality yields
    \[
      V_p^T U_k =
      \begin{bmatrix}
        V_k^T \\ \bar{V}_{p-k}^T
      \end{bmatrix}
      U_k =
      \begin{bmatrix}
        I_k \\ 0
      \end{bmatrix},
      \quad \text{and} \quad
      V_p V_p^T U_k =
      V_k.
    \]
    As in the first part of the proof, \(y\supL_k = T_{k-1,k}^T q\) for some \(q \in \R^{k-1}\), and therefore, \(x\supL_k = V_p V_p^T A^T U_{k-1} q\).
    Dual feasibility imposes that
    \begin{align*}
       0 \phantom{-} & = \phantom{-}U_{k-1}^T A V_p V_p^T U_p U_p^T (x\supL_k - x_\star)
      \\ & = \phantom{-}U_{k-1}^T A V_p U_p^T V_p \left( \begin{bmatrix} y\supL_k \\ 0 \end{bmatrix} - y_\star \right)
      \\ & = \phantom{-}U_{k-1}^T A (x\supL_k - x_\star)
      \\ & = -U_{k-1}^T r\supL_k,
    \end{align*}
    where we used biorthogonality, and~\eqref{eq:xLk-Vp}, and is satisfied because of~\eqref{eq:rkL}.
  \end{proof}

  Note that~\eqref{eq:xL-min-norm} continues to hold if the objective is measured in the \(U_p U_p^T\)-norm.
  Although this norm is no longer iteration dependent, it is unknown until the end of the biorthogonalization process.

  In the symmetric case, where \(V_k = U_k\) is orthogonal and \(T_k = T_k^T\), the \SYMMLQ iterate solves the problem
  \begin{equation}
    \label{eq:symmlq-min-norm}
    \minimize{x} \, \|x - x_\star\| \ \st \, x \in \Range(A V_{k-1}),
  \end{equation}
  which coincides with~\eqref{eq:xL-min-error}.

  \subsection{Numerical experiments}

\begin{figure}[ht]
    \centering
    \includetikzgraphics[width=0.8\textwidth]{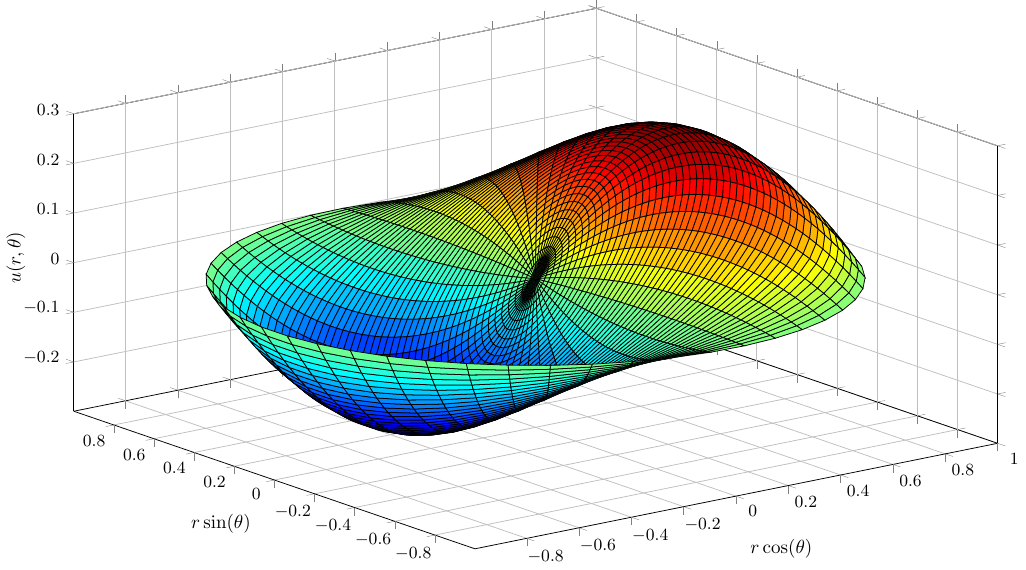}
    \label{fig:solution_u}
    \caption{%
      Solution $u(r, \theta) = r(1-r)\cos(\theta)$ of~\eqref{poisson} with $g(\theta) = 0$, $f(r, \theta) = - 3\cos(\theta)$ and $R=1$.
    }
  \end{figure}

  Non-homogeneous linear PDEs with variable coefficients of the form
  \begin{equation}
    \label{eq:linear_PDE}
    \displaystyle \sum_{i=1}^n \sum_{j=1}^p a_{i,j}(x) \frac{\partial^j u(x)}{\partial x_i^j} = b(x)
  \end{equation}
  are frequent when physical phenomena are modeled in polar, cylindrical or spherical coordinates.
  The discretization of~\eqref{eq:linear_PDE} often leads to a nonsymmetric square system.
  Such is the case with Poisson's equation $\Delta u = f$ used, for instance, to describe the gravitational or electrostatic field caused by a given mass density or charge distribution.
  The 2D Poisson equation in polar coordinates with Dirichlet boundary conditions is
  \begin{subequations}
    \label{poisson}
    \begin{alignat}{2}
      \frac{1}{r} \frac{\partial}{\partial r} \left( r \frac{\partial u(r, \theta)}{\partial r} \right) + \frac{1}{r^2} \frac{\partial^2 u(r, \theta)}{\partial \theta^2} & = f(r, \theta),
      \quad &
      (r,\theta) & \in (0, R) \times [0, 2\pi)
      \\
      u(R,\theta) & = g(\theta), \quad & \theta & \in [0, 2\pi),
    \end{alignat}
  \end{subequations}
  where \(R > 0\), the source term \(f\) and the boundary condition \(g\) are given.
  We discretize \eqref{poisson} using centered differences using \(50\) discretization points for \(r\) and \(50\) for \(\theta\), with $g(\theta) = 0$, $f(r, \theta) = - 3\cos(\theta)$ and $R=1$ so that~\eqref{poisson} models the response of an attached circular elastic membrane to a force.
  The resulting matrix has size \(2,500\) with $12,400$ nonzeros, and is block tridiagonal with extra diagonal blocks in the northeast and southwest corners.
  Each block on the main diagonal is tridiagonal but not symmetric.
  Each off-diagonal block is diagonal.
  More details on the discretization used are given by \cite{lai-2001}.
  The exact solution is represented in~\cref{fig:solution_u}.

\begin{figure}[ht]
    \centering
    \includetikzgraphics[width=0.49\textwidth]{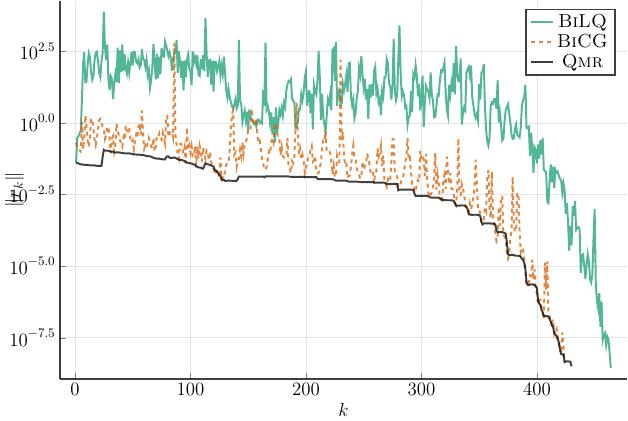}
    \hfill
    \includetikzgraphics[width=0.49\textwidth]{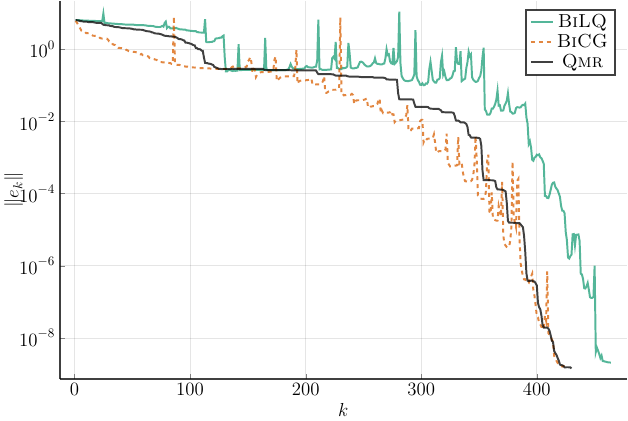}
    \label{fig:poisson}
    \caption{%
      Convergence curves of \BiLQ, \BiCG and \QMR iterates on~\eqref{poisson}.
      The figures show the residual (left) and error (right) history for each method.
    }
  \end{figure}

  We compare \BiLQ with our implementation of \QMR without lookahead.
  We also simulate \BiCG by way of the transition from \(x_k\supL\) to \(x_k\supC\) in \cref{alg:bilq}.
  \cref{fig:poisson} reports the residual and error history of \BiLQ, \BiCG and \QMR on~\eqref{poisson}.
  To compute $\|r_k\|$ and $\|e_k\|$, residuals $b - Ax_k$ and errors $x_k - x_{\star}$ are explicitly calculated at each iteration. We compute a reference solution with Julia's backslash command.
  We run each method with an absolute tolerance $\varepsilon_a = 10^{-10}$ and a relative tolerance $\varepsilon_r = 10^{-7}$ such that algorithms stop when $\|r_k\| \leq \varepsilon_a + \|b\| \varepsilon_r$.

  We also compare \BiLQ with \BiCG and \QMR on matrices SHERMAN5 and RAEFSKY1, with their respective right-hand side, from the UFL collection of \cite{davis-hu-2011}.\footnote{Now the SuiteSparse Matrix Collection \https{sparse.tamu.edu}.}
  System SHERMAN5 has size $3,312$ with $20,793$ nonzeros and RAEFSKY1 has size $3,242$ with $293,409$ nonzeros.
  A Jacobi preconditioner is used for both systems.

  \begin{figure}[ht]
    \centering
    \includetikzgraphics[width=0.49\textwidth]{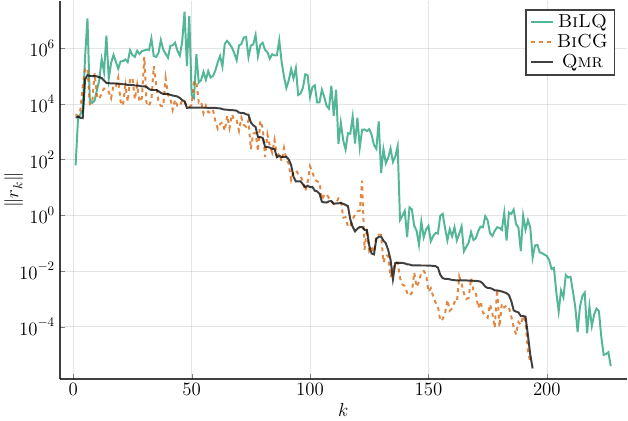}
    \hfill
    \includetikzgraphics[width=0.49\textwidth]{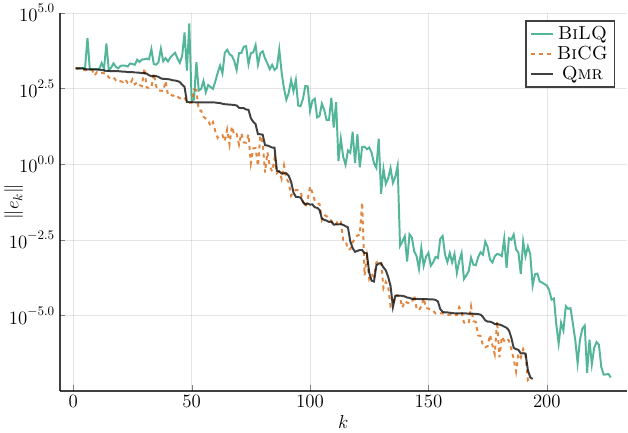}
    \label{fig:sherman5}
    \caption{%
      Convergence curves of \BiLQ, \BiCG and \QMR iterates for the SHERMAN5 system.
      The figures show the residual (left) and error (right) history for each method.
    }
  \end{figure}

  \begin{figure}[ht]
    \centering
    \includetikzgraphics[width=0.49\textwidth]{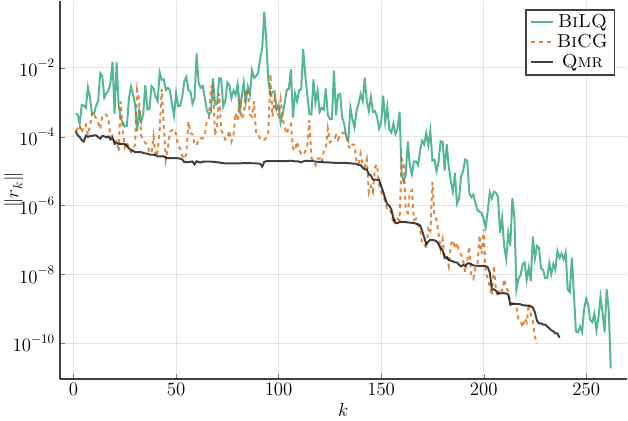}
    \hfill
    \includetikzgraphics[width=0.49\textwidth]{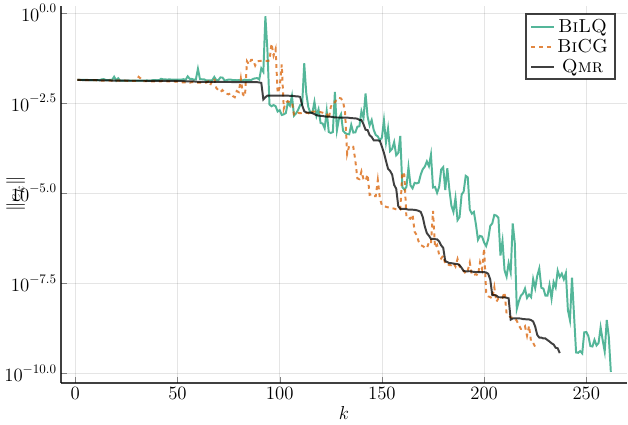}
    \label{fig:raefsky1}
    \caption{%
      Convergence curves of \BiLQ, \BiCG and \QMR iterates for the RAEFSKY1 system.
      The figures show the residual (left) and error (right) history for each method.
    }
  \end{figure}

  \cref{fig:poisson}, \cref{fig:sherman5} and \cref{fig:raefsky1} all show that in \BiLQ, neither the residual nor the error are monotonic in general.
  They also appear more erratic than those of \QMR.
  As in the symmetric case, both generally lag compared to those of \BiCG and \QMR, but are not far behind.
  We experimented with other systems and observed the same qualitative behavior.
  As showed in \cref{sec:bilq-properties}, although \BiLQ is a minimum-error-type method, this error is minimized over a different space than that where \(x\supL_k\) and \(x\supC_k\) reside---see \cref{prop:xL-min-norm}.
  This situation is analogous to that between \SYMMLQ and \CGM in the symmetric case \citep{estrin-orban-saunders-2019c}.
  Thus, the possibility of transferring to the \BiCG point, when it exists, is attractive.
  Because the \BiCG residual is easily computable, transferring based on the residual norm is readily implemented.
  The determination of upper bounds on the error suitable as stopping criteria remains the subject of active research \citep{estrin-orban-saunders-2019a,estrin-orban-saunders-2019b,estrin-orban-saunders-2019c}.

  \subsection{Discussion}

  Like \QMR, the \BiLQ iterate is well defined at each step even if $T_k$ is singular, whereas $x_k\supC$ is undefined when $\bar{\delta}_k= 0$.
  A simple example is
  \[
    A =
    \begin{bmatrix}
      0 & -1 \\
      1 & \phantom{-}1
    \end{bmatrix},
    \qquad
    b = c =
    \begin{bmatrix}
      1 \\
      0
    \end{bmatrix}.
  \]
  According to \cref{alg:bilanczos}, $\beta_1 = \gamma_1 = 1$, $v_1 = u_1 = b = c$.
  Then $\alpha_1 = u_1^T A v_1 = 0$, $T_1 = \begin{bmatrix} \alpha_1 \end{bmatrix}$ is singular, and $T_1 y_1 = \beta_1$ is inconsistent.
  \BiCG and its variants \CGS \citep{sonneveld-1989} and \BiCGSTAB \citep{van_der_vost-1992} all fail.
  However, $T_2$ is not singular and the \BiCG point exists, although we cannot compute it without lookahead. In finite precision arithmetic, such exact breakdown are rather rare. But near-breakdowns ($\bar{\delta}_k \approx 0$) may happen and lead to numerical instabilities in ensuing iterations.
  An additional drawback of \BiCG is that the LU decomposition of $T_k$ might not exist without pivoting even if $T_k$ is nonsingular whereas the LQ factorization of $T_{k-1,k}$ is always well defined.

  \section{Adjoint systems}
  \label{sec:adjoint-systems}

  Motivated by fuild dynamics applications, \cite{pierce-giles-2000} describe a method for doubling the order of accuracy of estimates of integral functionals involving the solution of a PDE.
  Consider a well-posed linear PDE \(L u = f\) on a domain \(\Omega\) subject to homogeneous boundary conditions, where \(L\) is a differential operator of the form~\eqref{eq:linear_PDE} and \(f \in L_2(\Omega)\).
  Suppose we wish to evaluate the functional \(J(u) := \langle u, \, g \rangle\), where \(g \in L_2(\Omega)\) and \(\langle \cdot, \, \cdot \rangle\) represents an integral inner product on \(L_2(\Omega)\).
  The problem may be stated equivalently as evaluating the functional \(\langle v, \, f \rangle\) where \(v\) solves the adjoint PDE \(L^* v = g\) because $\langle v, \, f \rangle = \langle v, \, Lu \rangle = \langle L^*v, \, u \rangle = \langle g, \, u \rangle$.
  
  Let the discretization of $L$ yield the linear system $A u_D = f_D$ with $D$ a set of points that define a grid on $\Omega$.
  For certain types of PDEs and certain discretization schemes, \(A^T\) is an appropriate discretization of \(L^*\).
  \cite{pierce-giles-2000} provide examples with linear operators such as Poisson's equation discretized by finite differences in 1D and by finite elements in 2D, but their discretizations are symmetric.
  Their method also applies to cases where \(A \neq A^T\) but in such cases, the discretization of the primal and dual equations commonly differ.
  Therefore, there is a need for methods that solve an unsymmetric primal system and its adjoint simultaneously.
  \cite{lu-darmofal-2003} and \cite{golub-stoll-wathen-2008} were also interested in this problem for scattering amplitude evaluation.
  \cite{lu-darmofal-2003} devise a modification of \QMR in which the two initial vectors are \(b\) and \(c\) and a quasi residual is minimized for both the primal and adjoint systems via an updated QR factorization.
  \cite{golub-stoll-wathen-2008} apply \USYMQR \citep{saunders-simon-yip-1988} to both the primal and the adjoint system\footnote{Although they call \USYMQR the ``generalized \LSQR''.} simultaneously by updating two QR factorizations.
  The advantage of their approach is that it produces monotonic residuals for both systems.

  Assume we use a method to compute \(u_D\) and to solve $A^T v_D = g_D$ such that $\|u - u_D\| \in O(h^p)$ and $\|v - v_D\| \in O(h^p)$, where \(h\) describes the grid coarseness.
  From $u_D$ and $v_D$ we compute approximations $u_h \approx u$ and $v_h \approx v$ over $\Omega$ by way of an interpolation of higher order than the discretization.
  Define $f_h := L u_h$ and $g_h := L^* v_h$.
  Instead of \(J(u) \approx \langle u_h, \, g \rangle\), an approximation of order \(p\), we may obtain one of order \(2p\) via the identity
  \begin{equation}
    \label{eq:superconvergent-approx}
    \langle g, \, u \rangle =
    \langle g, \, u_h \rangle -
    \langle v_h, \, f_h - f \rangle +
    \langle g_h - g, \, u_h - u \rangle.
  \end{equation}
  The first two terms constitute our new approximation while the remaining error term can be expressed as $\langle g_h - g, L^{-1}(f_h - f) \rangle = O(h^{2p})$.

  From this point, we consider, in addition to~\eqref{eq:Ax=b}, the adjoint system
  \begin{equation}
    \label{eq:adjoint-system}
    A^T t = c.
  \end{equation}
  Solving simultaneously primal and dual systems can also be formulated as solving the symmetric and indefinite system
  \begin{equation}
    \label{eq:adjoint_systems}
    \begin{bmatrix}
      0   & A \\
      A^T & 0
    \end{bmatrix}
    \begin{bmatrix}
    t \\
    x
    \end{bmatrix}
    =
    \begin{bmatrix}
    b \\
    c
    \end{bmatrix}.
  \end{equation}
  \MINRES or \MINRESQLP \citep{choi-paige-minres-2011} are prime candidates for~\eqref{eq:adjoint_systems} and will serve as a basis for comparison.

  In the context of \cref{alg:bilanczos}, we can take advantage of the two initial vectors \(b\) and \(c\) to combine \BiLQ and \QMR and solve both the primal and adjoint systems simultaneously at no other extra cost than that of updating solution and residual estimates.
  We call the resulting method \BiLQR.
  Contrary to the approach of \cite{lu-darmofal-2003}, no extra factorization updates are necessary.
  Instead of approximating \(u_D\) and \(v_D\) by minimizing two quasi residuals, \BiLQR minimizes one quasi residual and computes the second approximation via a minimum-norm subproblem.

  A similar method based on the orthogonal tridiagonalization process of \cite{saunders-simon-yip-1988} can be derived by combining \USYMLQ and \USYMQR, which we call \TriLQR, and which is to the approach of \cite{golub-stoll-wathen-2008} as \BiLQR is to that of \cite{lu-darmofal-2003}.
  \TriLQR remains well defined for rectangular \(A\).

  \subsection{Description of \BiLQR}

  \BiLQR updates an approximate solution $t\supQ_{k-1} = U_{k-1} f\supQ_{k-1}$ of $A^T t = c$ by solving the \QMR least-squares subproblem
  \begin{equation}
    \label{eq:qmr-sub-dual}
    \minimize{f} \ \|T_{k-1,k}^T f - \gamma_1 e_1\|
    \quad \Longleftrightarrow \quad
    \minimize{f} \
    \left\|
      \begin{bmatrix}
        L_{k-1}^T \\
        0
      \end{bmatrix} f -
      Q_k \gamma_1 e_1
    \right\|
  \end{equation}
  because the QR factorization of \(T_{k-1,k}^T\) is readily available.
  Define $\bar{h}_k = Q_k \gamma_1 e_1 = (h_{k-1}, \, \bar{\psi}_k) = (\psi_1, \, \cdots, \, \psi_{k-1}, \, \bar{\psi}_k)$.
  The components of $\bar{h}_k$ are updated according to
  \begin{subequations}
    \label{hk}
    \begin{align}
      \bar{\psi}_1 & = \gamma_1, \\
      \psi_k & = c_{k+1} \bar{\psi}_{k},~k \ge 1, \\
      \bar{\psi}_{k+1} & = s_{k+1} \bar{\psi}_{k}, \quad k \ge 1.
    \end{align}
  \end{subequations}
  The solution of~\eqref{eq:qmr-sub-dual} is $f\supQ_{k-1} = L_{k-1}^{-T} h_{k-1}$ and the least-squares residual norm is $|\bar{\psi}_k|$.
  To avoid storing $U_k$, we define $W_k = U_k L_k^{-T}$, which can be updated as
  \begin{subequations}
    \label{Wk}
    \begin{align}
      w_1 & = u_1 / \delta_1, \\
      w_2 & = (u_2 - \lambda_1 w_1) / \delta_2, \\
      w_k & = (u_k - \lambda_{k-1} w_{k-1} - \varepsilon_{k-2} w_{k-2}) / \delta_k, \quad k \ge 3.
    \end{align}
  \end{subequations}
  At the next iteration, $t\supQ_k$ can be recursively updated according to
  \[
    t\supQ_k =
    U_k f\supQ_k =
    U_k L_k^{-T} h_k =
    W_k h_k = W_{k-1} h_{k-1} + \psi_k w_k = t\supQ_{k-1} + \psi_k w_k.
  \]

  The \QMR residual is
  \[
    r_k\supQ =
    c - A^T t\supQ_k =
    U_{k+1}(\gamma_1 e_1 - T_{k,k+1}^T f\supQ_k) =
    \bar{\psi}_{k+1} U_{k+1} Q_{k+1}^T e_{k+1}^T,
  \]
  so that
  \[
    \|r_k\supQ\| \leq
    \|U_{k+1}\|_F \, \|\bar{\psi}_{k+1} Q_{k+1}^T e_{k+1}^T\| \leq
    \|\bar{\psi}_{k+1}\|\sqrt{\tau_{k+1}},
  \]
  where \(\tau_{k+1} = \sum_{i=1}^{k+1} \|u_i\|^2 = \tau_k + \|u_{k+1}\|^2\).
  If the $u_k$ are normalized, then $\tau_k = k$.
  \autoref{alg:bilqr} states the complete procedure.

  The following result states a minimization property of the \QMR residual in an iteration-dependent norm.

  \begin{shadyproposition}
    \label{prop:xQ-min-res}
    The \((k-1)\)th \QMR iterate \(t\supQ_{k-1}\) solves
    \begin{equation}
      \label{eq:xQ-min-res}
      \minimize{t} \, \|c - A^T t\|_{V_k V_k^T} \ \st \, t \in \Range(U_{k-1}).
  \end{equation}
  In addition, \(\|r\supQ_k\|_{V_k V_k^T}\) is monotonically decreasing.
  \end{shadyproposition}

  \begin{proof}
    The set of constraints of~\eqref{eq:xQ-min-res} imposes that there exist \(f \in \R^{k-1}\) such that \(t = U_{k-1} f\).
    By biorthogonality, the objective value at such an \(t\) can be written \(\|c- A^T U_{k-1} f\|_{V_k V_k^T} = \|c - U_k T_{k-1,k}^T f\|_{V_k V_k^T} = \|\gamma_1 e_1 - T_{k-1,k}^T f\|\). We recover the subproblem~\eqref{eq:qmr-sub-dual}.

  For the second part, $\|r\supQ_k\|_{V_{k+1} V_{k+1}^T} = |\bar{\psi}_{k+1}| = |s_{k+1}| |\bar{\psi}_k| = |s_{k+1}|\|r\supQ_{k-1}\|_{V_k V_k^T}$.
  \end{proof}

  Note that \cref{prop:xQ-min-res} continues to hold if \(r\supQ_k\) is measured in the \(V_p V_p^T\)-norm.

  \subsection{Description of \TriLQR}

  The \cite{saunders-simon-yip-1988} tridiagonalization process generates sequences of vectors $\{v_k\}$ and $\{u_k\}$ such that $v_i^T v_j = \delta_{ij}$ and $u_i^T u_j = \delta_{ij}$ in exact arithmetic for as long as the process does not break down. The process is summarized as \autoref{alg:ssy_process}.

  \begin{algorithm}[ht]
    \caption{%
      Tridiagonalization Process
    }
    \label{alg:ssy_process}
    \begin{algorithmic}[1]
      \Require $A$, $b$, $c$
      \State $v_0 = 0$, $u_0 = 0$
      \State $\beta_1 v_1 = b$, $\gamma_1 u_1 = c$ \Comment{$(\beta_1,~\gamma_1) > 0$ so that $\|v_1\| = \|u_1\| = 1$}
      \For{$k$ = 1,~2,~\(\dots\)}
        \State $q = A u_k - \gamma_k v_{k-1}$, $\alpha_k = v_k^T q$
        \State $p = A^T v_k - \beta_k u_{k-1}$
        \State $\beta_{k+1} v_{k+1} = q - \alpha_k v_k$ \Comment{$\beta_{k+1} > 0$ so that $\|v_{k+1}\| = 1$}
        \State $\gamma_{k+1} u_{k+1} = p - \alpha_k u_k$ \Comment{$\gamma_{k+1} > 0$ so that $\|u_{k+1}\| = 1$}
      \EndFor
    \end{algorithmic}
  \end{algorithm}

  At the end of the $k$-th iteration, we have
  \begin{subequations}
    \label{eq:ssy_process}
    \begin{align}
         A U_k & = V_k T_k~ + \beta_{k+1} v_{k+1} e_k^T = V_{k+1} T_{k+1,k}
      \\ A^T V_k & = U_k T_k^T + \gamma_{k+1} u_{k+1} e_k^T = U_{k+1} T_{k,k+1}^T,
    \end{align}
  \end{subequations}
  to be compared with~\eqref{eq:bilanczos}.

  \cite{saunders-simon-yip-1988} develop two methods based on \autoref{alg:ssy_process}.
  \USYMLQ generates an approximation to a solution of~\eqref{eq:Ax=b} of the form \(x\supLQ_k = U_k y\supLQ_k\), where \(y\supLQ_k \in \R^k\) solves
  \begin{equation}
    \label{eq:usymlq-sub}
    \minimize{y} \ \|y\| \quad \st \ T_{k-1,k} y = \beta_1 e_1.
  \end{equation}
  With~\eqref{eq:ssy_process} and~\eqref{eq:usymlq-sub}, we have the following analogue of \cref{prop:xL-min-norm} and~\eqref{eq:symmlq-min-norm}.

  \begin{shadyproposition}
    \label{prop:xLQ-min-norm}
    Let \(x_\star\) be the exact solution of~\eqref{eq:Ax=b}.
    The \(k\)th \USYMLQ iterate \(x\supLQ_k\) solves
    \begin{equation}
      \label{eq:xLQ-min-norm}
      \minimize{x} \, \|x\| \ \st \, x \in \Range(U_k), \ b - A x \perp \Range(U_{k-1}),
    \end{equation}
    and
    \begin{equation}
      \label{eq:xLQ-min-error}
      \minimize{x} \, \|x - x_\star\| \ \st \, x \in \Range(A^T V_{k-1}).
    \end{equation}
  \end{shadyproposition}

  \begin{proof}
    The proof is nearly identical to that of \cref{prop:xL-min-norm} and relies on the fact that \(r\supLQ_k := b - A x\supLQ_k\) is a combination of \(u_k\) and \(u_{k+1}\) \cite[\S3.2.2]{buttari-orban-ruiz-titley_peloquin-2019}.
  \end{proof}

  The second method, \USYMQR, generates an approximation \(t\supQR_k = V_k f\supQR_k\) where \(f\supQR_k \in \R^k\) solves
  \begin{equation}
    \label{eq:usymqr-sub}
    \minimize{f} \ \|T_{k,k+1}^T f - \gamma_1 e_1\|.
  \end{equation}
  The following property applies to \(t\supQR_k\) due to our assumption that~\eqref{eq:Ax=b} is consistent.

  \begin{shadyproposition}[{\protect \citealp[Theorem~$1$]{buttari-orban-ruiz-titley_peloquin-2019}}]
    \label{prop:tQR-min-LS}
    Assume \(b \in \Range(A)\).
    Then \USYMQR finds the minimum-norm solution of
    \[
      \minimize{t} \ \|A^T t - c\|.
    \]
  \end{shadyproposition}

  Of course, \(A\) nonsingular implies that the solution to~\eqref{eq:adjoint-system} is unique but \cref{prop:tQR-min-LS} applies more generally to rectangular and/or rank-deficient \(A\). 

  When $A = A^T$ and $b = c$, \cref{alg:ssy_process} coincides with the symmetric Lanczos process, and \USYMLQ and \USYMQR are equivalent to \SYMMLQ and \MINRES \citep{paige-saunders-1975}, respectively.
  Besides the orthogonalization process, differences between those methods and \BiLQ and \QMR are the definition of $\bar{D}_k$ and $W_k$, and the fact that $u_k$ and $v_k$ are swapped.
  If stopping criteria are based on residual norms, expressions derived for methods based on \cref{alg:bilanczos} apply to methods based on \cref{alg:ssy_process}, but their expressions can simplified because $V_k$ and $U_k$ are orthogonal.
  \USYMQR and \USYMLQ can be combined into \TriLQR to solve both the primal and ajoint system simultaneously.
  We summarize the complete procedure as \autoref{algo:trilqr} and highlight lines with differences between the two algorithms.

  \begin{doublealgorithm}
  \begin{minipage}[t]{0.49\textwidth}
    \setlength{\intextsep}{0pt}
    \begin{algorithm}[H]
      \caption{\BiLQR}
      \label{alg:bilqr}
      \begin{algorithmic}[1]
      \Require $A$, $b$, $c$
      \State $\beta_1 v_1 = b$, $\gamma_1 u_1 = c$
      \State $\alpha_1 = u_1^T A v_1$
      \State \highlightline{$\beta_2 v_2 = A v_1 - \alpha_1 v_1$}
      \State \highlightline{$\gamma_2 u_2 = A^T u_1 - \alpha_1 u_1$}
      \State $c_1 = -1$, $s_1 = 0$, $\bar{\delta}_1 = \alpha_1$
      \State $\eta_1 = \beta_1$, \highlight{$\bar{d}_1 = v_1$}, $\bar{\psi}_1 = \gamma_1$
      \State $x\supL_1 = 0$, $t\supQ_0 = 0$
      \For{$k$ = 2,~3,~\(\dots\)}
        \State \highlightlinefor{$q = A v_k - \gamma_k v_{k-1}$, $\alpha_k = u_k^T q$}
        \State \highlightlinefor{$p = A^T u_k - \beta_k u_{k-1}$}
        \State $\beta_{k+1} v_{k+1} = q  - \alpha_k v_k$
        \State $\gamma_{k+1} u_{k+1} = p - \alpha_k u_k$
        \State $\delta_{k-1} = (\bar{\delta}_{k-1}^2 + \gamma_k^2)^{\frac12}$
        \State $c_k = \bar{\delta}_{k-1} / \delta_{k-1}$
        \State $s_k = \gamma_k / \delta_{k-1}$
        \State $\varepsilon_{k-2} = s_{k-1} \beta_k$
        \State $\lambda_{k-1} = - c_{k-1} c_k \beta_k + s_k \alpha_k$
        \State $\bar{\delta}_k = - c_{k-1} s_k \beta_k - c_k \alpha_k$
        \State $\zeta_{k-1} = \eta_{k-1} / \delta_{k-1}$
        \State $\eta_k = -\varepsilon_{k-2} \zeta_{k-2} - \lambda_{k-1} \zeta_{k-1}$
        \State \highlightlinefor{$d_{k-1} = c_k \bar{d}_{k-1} + s_k v_k$}
        \State \highlightlinefor{$\bar{d}_k = s_k \bar{d}_{k-1} - c_k v_k$}
        \State $\psi_{k-1} = c_k \bar{\psi}_{k-1}$
        \State $\bar{\psi}_k = s_k \bar{\psi}_{k-1}$
        \State \highlightlinefor{$w_{k-1} = \frac{u_{k-1} - \lambda_{k-2} w_{k-2} - \varepsilon_{k-3} w_{k-3}}{\delta_{k-1}}$}
        \State $x\supL_k = x\supL_{k-1} + \zeta_{k-1} d_{k-1}$
        \State $t\supQ_{k-1} = t\supQ_{k-2} + \psi_{k-1} w_{k-1}$
      \EndFor
      \If{$\bar{\delta}_k \neq 0$}
        \State $\bar{\zeta}_k = \eta_k / \bar{\delta}_k$
        \State $x\supC_k = x\supL_k + \bar{\zeta}_k \bar{d_k}$
      \EndIf
      \end{algorithmic}
    \end{algorithm}
  \end{minipage}
  \hfill
  \begin{minipage}[t]{0.49\textwidth}
    \setlength{\intextsep}{0pt}
    \begin{algorithm}[H]
    \caption{\TriLQR}
    \label{algo:trilqr}
    \begin{algorithmic}
      \Require $A$, $b$, $c$
      \State $\beta_1 v_1 = b$, $\gamma_1 u_1 = c$
      \State $\alpha_1 = u_1^T A v_1$
      \State \highlightline{$\beta_2 v_2 = A   u_1 - \alpha_1 v_1$}
      \State \highlightline{$\gamma_2 u_2 = A^T v_1 - \alpha_1 u_1$}
      \State $c_1 = -1$, $s_1 = 0$, $\bar{\delta}_1 = \alpha_1$
      \State $\bar{\eta}_1 = \beta_1$, \highlight{$\bar{d}_1 = u_1$}, $\bar{\psi}_1 = \gamma_1$
      \State $x\supLQ_1 = 0$, $t\supQR_0 = 0$
      \For{$k$ = 2,~3,~\(\dots\)}
        \State \highlightlinefor{$q = A u_k - \gamma_k v_{k-1}$, $\alpha_k = v_k^T q$}
        \State \highlightlinefor{$p = A^T v_k - \beta_k u_{k-1}$}
        \State $\beta_{k+1} v_{k+1} = q  - \alpha_k v_k$
        \State $\gamma_{k+1} u_{k+1} = p - \alpha_k u_k$
        \State $\delta_{k-1} = (\bar{\delta}_{k-1}^2 + \gamma_k^2)^{\frac12}$
        \State $c_k = \bar{\delta}_{k-1} / \delta_{k-1}$
        \State $s_k = \gamma_k / \delta_{k-1}$
        \State $\varepsilon_{k-2} = s_{k-1} \beta_k$
        \State $\lambda_{k-1} = - c_{k-1} c_k \beta_k + s_k \alpha_k$
        \State $\bar{\delta}_k = - c_{k-1} s_k \beta_k - c_k \alpha_k$
        \State $\zeta_{k-1} = \eta_{k-1} / \delta_{k-1}$
        \State $\eta_k = -\varepsilon_{k-2} \zeta_{k-2} - \lambda_{k-1} \zeta_{k-1}$
        \State \highlightlinefor{$d_{k-1} = c_k \bar{d}_{k-1} + s_k u_k$}
        \State \highlightlinefor{$\bar{d}_k = s_k \bar{d}_{k-1} - c_k u_k$}
        \State $\psi_{k-1} = c_k \bar{\psi}_{k-1}$
        \State $\bar{\psi}_k = s_k \bar{\psi}_{k-1}$
        \State \highlightlinefor{$w_{k-1} = \frac{v_{k-1} - \lambda_{k-2} w_{k-2} - \varepsilon_{k-3} w_{k-3}}{\delta_{k-1}}$}
        \State $x\supLQ_k = x\supLQ_{k-1} + \zeta_{k-1} d_{k-1}$
        \State $t\supQR_{k-1} = t\supQR_{k-2} + \psi_{k-1} w_{k-1}$
      \EndFor
      \If{$\bar{\delta}_k \neq 0$}
        \State $\bar{\zeta}_k = \eta_k / \bar{\delta}_k$
        \State $x\supCG_k = x\supLQ_k + \bar{\zeta}_k \bar{d_k}$
      \EndIf
    \end{algorithmic}
  \end{algorithm}
  \end{minipage}
  \end{doublealgorithm}

  \BiLQR and \TriLQR both need nine $n$-vectors: $u_k$, $u_{k-1}$, $v_k$, $v_{k-1}$, $w_k$, $w_{k-1}$, $\bar{d}_k$, $x_k$ and $t_{k-1}$ whereas \MINRESQLP applied to \eqref{eq:adjoint_systems} can be implemented with five $(2n)$-vectors.
  Two more $n$-vectors are needed when in-place ``gemv'' updates are not explicitly available.
Table~\ref{tab:cost_bilqr_trilqr_minresqlp} summarizes the cost of \BiLQR, \TriLQR, \MINRESQLP and variants from \cite{lu-darmofal-2003} and \cite{golub-stoll-wathen-2008}, developed for adjoint systems.
  An advantage of \MINRESQLP and \TriLQR is that adjoint systems can be solved even if $b^T c = 0$, which is not possible with \BiLQR.
  In addition, serious breakdowns $q^T p = 0$ with $p \neq 0$ and $q \neq 0$ are not a problem with \TriLQR.
  \TriLQR is similar in spirit to the recent method \USYMLQR of \cite{buttari-orban-ruiz-titley_peloquin-2019} for solving symmetric saddle-point systems, but is slightly cheaper.

  \begin{table}[ht]
    \label{tab:cost_bilqr_trilqr_minresqlp}
    \caption{Storage and cost per iteration of methods for solving~\eqref{eq:Ax=b} and~\eqref{eq:adjoint-system} simultaneously.}
    \begin{center}
      \setlength{\arrayrulewidth}{1pt}
      \begin{tabular}{lrrrr}
        \hline
        & $n$-vectors & dots & scal & axpy \\
        \hline
        \BiLQR                          & 9  & 2   & 5   & 10 \\
        \TriLQR                         & 9  & 2   & 5   & 10 \\
        \MINRESQLP                      & 10 & 4   & 8   & 14 \\
        \cite{lu-darmofal-2003}         & 10 & 2   & 6   & 10 \\
        \cite{golub-stoll-wathen-2008}  & 10 & 2   & 6   & 10 \\
        \hline
      \end{tabular}
    \end{center}
  \end{table}

  \subsection{Applications}

  For the purpose of a simple illustration, we consider a one-dimensional ODE and a two-dimensional PDE.
  Consider first the linear ODE with constant coefficients
  \begin{subequations}
    \label{eq:primal-1D}
    \begin{alignat}{2}
      \chi_1 u''(x) + \chi_2 u'(x) + \chi_3 u(x) &= f(x)
      \quad & x & \in \Omega \\
      u(x) &= 0 \quad & x & \in \partial\Omega,
    \end{alignat}
  \end{subequations}
  where \(\Omega = [0, \, 1]\), and say we are interested in the value of the linear functional
  \begin{equation}
    \label{eq:primal-functional}
    J(u) = \int_\Omega u(x) g(x) \, \mathrm{d}\Omega,
  \end{equation}
  where \(u\) solves~\eqref{eq:primal-1D} and \(g \in L_2(\Omega)\).
  The adjoint equation can be derived from~\eqref{eq:primal-1D} using integration by parts:
  \begin{subequations}
    \label{eq:dual-1D}
    \begin{alignat}{2}
      \chi_1 v''(x) - \chi_2 v'(x) + \chi_3 v(x) &= g(x)
      \quad & x & \in \Omega \\
      v(x) &= 0 \quad & x & \in \partial\Omega.
    \end{alignat}
  \end{subequations}
  Note that the only difference between the primal and adjoint equations resides in the sign of odd-degree derivatives.
  The discussion in \cref{sec:adjoint-systems} ensures that
  \begin{equation}
    \label{eq:dual-functional}
    G(v) := \int_\Omega f(x) v(x) \, \mathrm{d}\Omega = J(u).
  \end{equation}

  Consider the uniform discretization $x_i = ih$, $i = 0, \dots, N+1$, where $h = 1 / (N+1)$.
  We use centered finite differences of order 2, i.e.,
  \[
    u'(x_i) = \frac{u_{i+1} - u_{i-1}}{2h} + O(h^2),
    \quad
    u''(x_i) = \frac{u_{i-1} - 2u_{i} + u_{i+1}}{h^2} + O(h^2).
  \]
  We obtain $u(x_i)$ for $x_i \in D := \{x_i \mid i \in 1,\dots,N\}$ from the tridiagonal linear system
  \[
    \begin{bmatrix}
      - 2\chi_1 + \chi_3 h^2  & \chi_1 + \chi_2 h      &        & \\
      \chi_1 - \chi_2 h       & - 2\chi_1 + \chi_3 h^2 & \ddots & \\
                                                       & \ddots   & \ddots & \chi_1 + \chi_2 h \\
                                                       &          & \chi_1 - \chi_2 h & - 2\chi_1 + \chi_3 h^2
    \end{bmatrix}
    \begin{bmatrix}
      u(x_1) \\
      \vdots \\
      \vdots \\
      u(x_N)
    \end{bmatrix}
    =
    h^2
    \begin{bmatrix}
      f(x_1) \\
      \vdots \\
      \vdots \\
      f(x_N)
    \end{bmatrix}.
  \]
  More compactly, we write $A u_D = f_D$.
  Similarly, we compute $v(x_i)$ for $x_i \in D$ from $A^T v_D = g_D$.
  Next, we compute an approximation of $u$ and $v$ over $\Omega$ by cubic spline interpolation, and the resulting functions are denoted $u_h$ and $v_h$.
  We impose that \(L u_h = f\) and \(L^* v_h = g\) on \(\partial\Omega\).
  We subsequently obtain \(f_h(x) := \chi_1 u_h''(x) + \chi_2 u_h'(x) + \chi_3 u_h(x)\).
  The end points conditions of the cubic splines impose that \(f_h\)
  coincide with \(f\)
  on \(\partial\Omega\).
  Finally, we compute the improved estimate~\eqref{eq:superconvergent-approx} using a three-point Gauss quadrature to approximate each
  \[
    \int_{x_i}^{x_{i+1}} g(x) u_h(x) \, \mathrm{d}x -
    \int_{x_i}^{x_{i+1}} v_h(x) (f_h(x) - f(x)) \, \mathrm{d}x
  \]
  on each subinterval to ensure that the numerical quadrature errors are smaller than the discretization error.

\begin{figure}[ht]
    \centering
    \includetikzgraphics[width=0.49\textwidth]{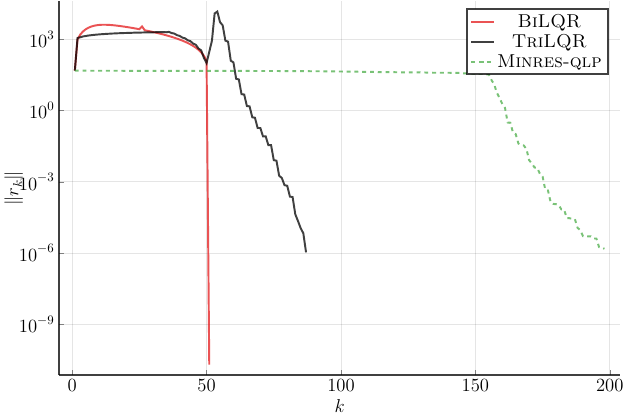}
    \hfill
    \includetikzgraphics[width=0.49\textwidth]{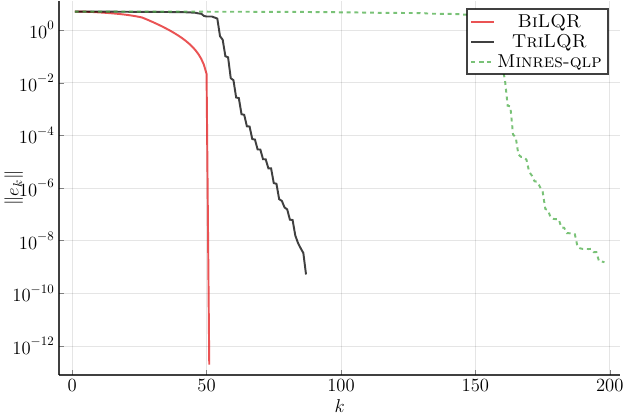}
    \label{fig:primal-ode}
    \caption{Residuals and errors norms of \BiLQR, \TriLQR and \MINRESQLP iterates for on~\eqref{eq:primal-1D}.}
  \end{figure}

\begin{figure}[ht]
    \centering
    \includetikzgraphics[width=0.49\textwidth]{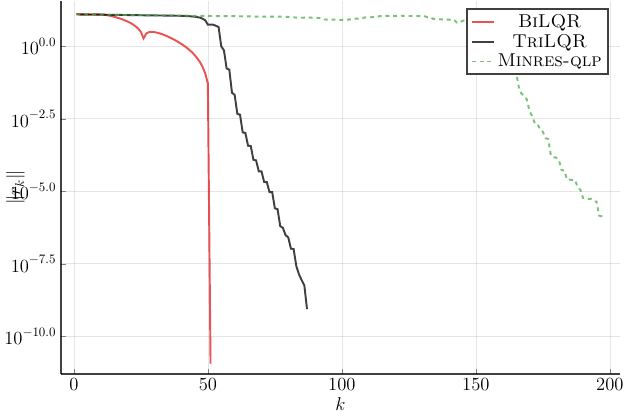}
    \hfill
    \includetikzgraphics[width=0.49\textwidth]{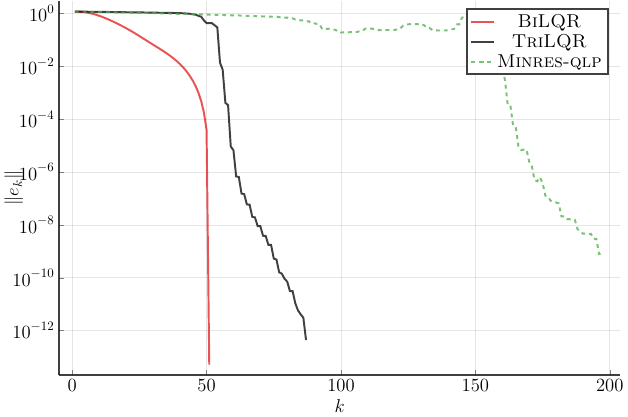}
    \label{fig:dual-ode}
    \caption{Residuals and errors norms of \BiLQR, \TriLQR and \MINRESQLP iterates on~\eqref{eq:dual-1D}.}
  \end{figure}

\begin{figure}[ht]
    \centering
    \includetikzgraphics[width=0.5\textwidth]{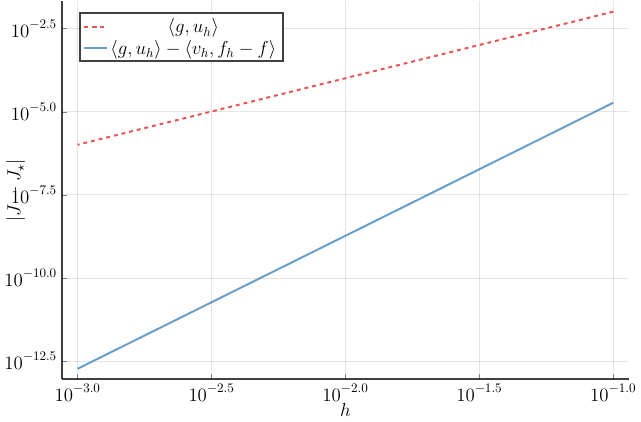}
    \hfill
    \includetikzgraphics[width=0.48\textwidth]{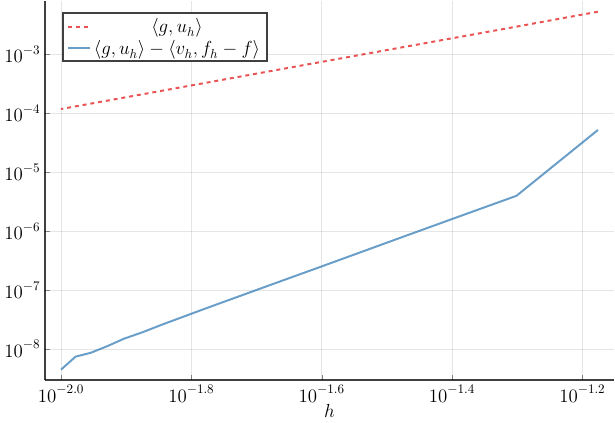}
    \label{fig:J}
    \caption{%
      Functional evaluation errors for~\eqref{eq:primal-1D}--\eqref{eq:dual-1D} (left) and~\eqref{primal-advection-diffusion}--\eqref{dual-advection-diffusion} (right).
    }
  \end{figure}

  We choose $n=50$, $\chi_1 = \chi_2 = \chi_3 = 1$, $g(x) = e^x$ and $f(x)$ such that the exact solution of~\eqref{eq:primal-1D} is $u_\star(x) = sin(\pi x)$. The resulting linear system has dimension $50$ with $148$ nonzeros. Those parameters ensure that $J_{\star} = \langle g, \, u_\star \rangle = (\pi (e + 1)) / (\pi^2 + 1)$. \cref{fig:primal-ode,fig:dual-ode} report the evolution of the residual and error on~\eqref{eq:Ax=b} and~\eqref{eq:adjoint-system} for~\eqref{eq:primal-1D} and~\eqref{eq:dual-1D}, respectively.
  \BiLQR terminates in 51 iterations, \TriLQR in 87 iterations and \MINRESQLP in 198 iterations.
  The left plot of \cref{fig:J} illustrates the error in the evaluation of \(J(u)\) as a function of \(h\) using the naive \(J(u) \approx J(u_h)\) and improved~\eqref{eq:superconvergent-approx} approximations.

  The steady-state convection-diffusion equation with constant coefficients
  \begin{subequations}
    \label{primal-advection-diffusion}
    \begin{alignat}{2}
      \kappa_1 \Delta u(x) + \kappa_2 \nabla \cdot u(x) & = f(x) \quad & x & \in \Omega
      \\
      u(x) & = 0 \quad & x & \in \partial \Omega,
    \end{alignat}
  \end{subequations}
  where \(f \in L_2(\Omega)\), describes the flow of heat, particles, or other physical quantities in situations where there is both diffusion and convection or advection.
  Assume as before that we are interested in the linear functional~\eqref{eq:primal-functional}.
  The adjoint equation of~\eqref{primal-advection-diffusion}, again obtained via integration by parts, reads
  \begin{subequations}
    \label{dual-advection-diffusion}
    \begin{alignat}{2}
      \kappa_1 \Delta v(x) - \kappa_2 \nabla \cdot v(x) & = g(x) \quad & x & \in \Omega
      \\
      v(x) & = 0 \quad & x & \in \partial \Omega,
    \end{alignat}
  \end{subequations}
  and duality ensures~\eqref{eq:dual-functional}. 

  In the case of heat transfer, $u(x)$ represents temperature and $f(x)$ sources or sinks.
  For example, with $g(x) = 1 / \vol(\Omega)$, $J(u)$ represents the average temperature in \(\Omega\).

  We choose $\Omega = [0, \, 1] \times [0, \, 1]$ and descretize~\eqref{primal-advection-diffusion} on a uniform $N \times N$ grid with the finite difference method such that the step along both coordinates is $h = 1 / (N+1)$.
  With centered second-order differences for first and second derivatives, the discretized operator has the structure
  \[
    A =
    \begin{bmatrix}
      T      & D_U    &        & \\
      D_L    & T      & \ddots & \\
             & \ddots & \ddots & D_U    \\
             &        & D_L    & T
    \end{bmatrix},
    \mathhfill
    T =
    \begin{bmatrix}
      \text{\small $-4 \kappa_1$} & \text{\small $\kappa_1 + \tfrac{1}{2} \kappa_2 h$} & & \\
      \text{\small $\kappa_1 - \tfrac{1}{2} \kappa_2 h$} & \text{\small $-4 \kappa_1$} & \ddots & \\
      & \ddots & \ddots & \text{\small $\kappa_1 + \tfrac{1}{2} \kappa_2 h$} \\
      & & \text{\small $\kappa_1 - \tfrac{1}{2} \kappa_2 h$} & \text{\small $-4 \kappa_1$}
    \end{bmatrix},
  \]
  $D_U = \diag(\kappa_1 + \tfrac{1}{2}\kappa_2 h)$, $D_L = \diag(\kappa_1 - \tfrac{1}{2} \kappa_2 h)$, where the right-hand sides \(b\) and \(c\) include the \(h^2\) term. Solutions $u_D$ and $v_D$ contain an approximation of $u$ and $v$ at grid points stored column by column.
  The discretization of~\eqref{dual-advection-diffusion} with the same scheme yields $A^T$. We compare \BiLQR, \TriLQR and \MINRESQLP on~\eqref{primal-advection-diffusion} and~\eqref{dual-advection-diffusion} with $\kappa_1=5$, $\kappa_2=20$, $N=50$, $g(x,y) = e^{x+y}$ and $f(x,y)$ such that the exact solution of~\eqref{primal-advection-diffusion} is $u_{\star}(x,y) = \sin(\pi x)\sin(\pi y)$. The resulting linear system has dimension $2,500$ with $12,300$ nonzeros.
  We use an absolute tolerance $\varepsilon_a = 10^{-10}$ and a relative tolerance $\varepsilon_r = 10^{-7}$, and terminate when both $\|r_k\| \leq \varepsilon_a + \|b\| \varepsilon_r$ for~\eqref{eq:Ax=b} and $\|r_k\| \leq \varepsilon_a + \|c\| \varepsilon_r$ for~\eqref{eq:adjoint-system} hold. 

\begin{figure}[ht]
    \centering
    \includetikzgraphics[width=0.49\textwidth]{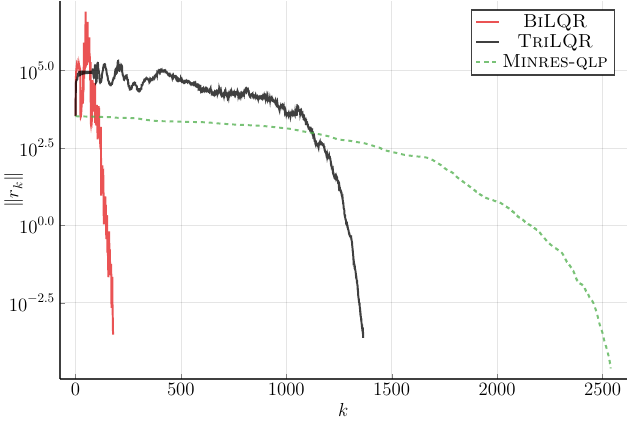}
    \hfill
    \includetikzgraphics[width=0.49\textwidth]{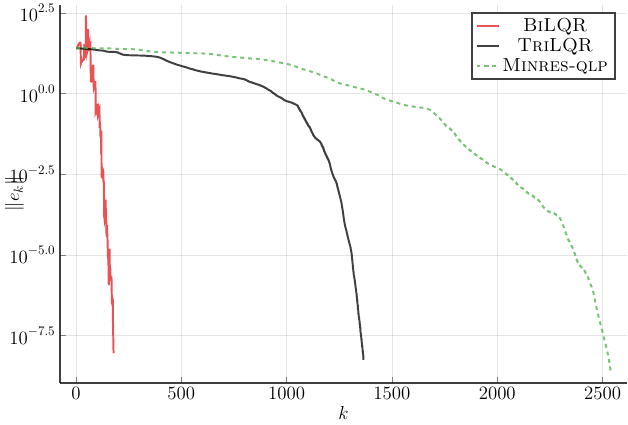}
    \label{fig:primal-advection-diffusion}
    \caption{Residuals and errors norms of \BiLQR, \TriLQR and \MINRESQLP iterates for on~\eqref{primal-advection-diffusion}.}
  \end{figure}

\begin{figure}[ht]
    \centering
    \includetikzgraphics[width=0.49\textwidth]{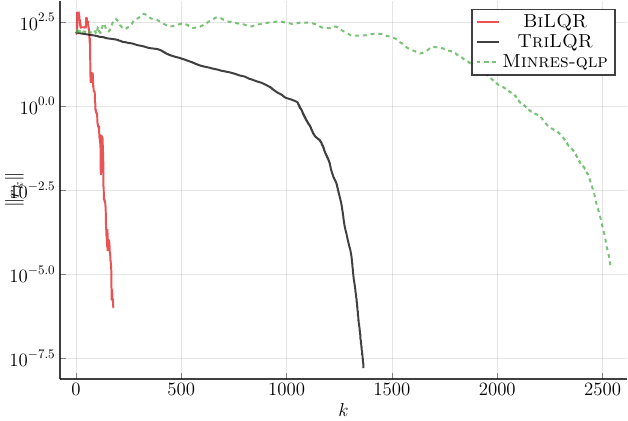}
    \hfill
    \includetikzgraphics[width=0.49\textwidth]{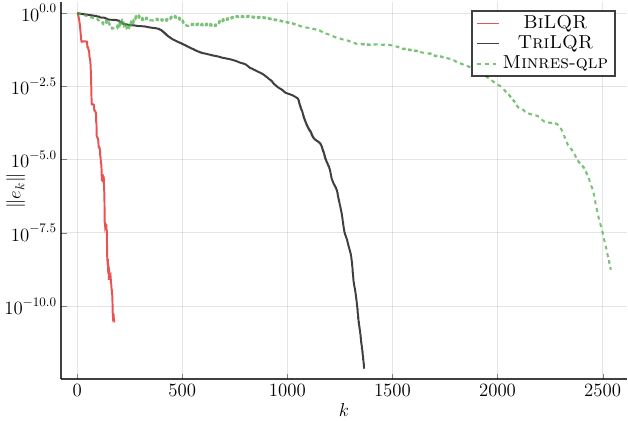}
    \label{fig:dual-advection-diffusion}
    \caption{Residuals and errors norms of \BiLQR, \TriLQR and \MINRESQLP iterates on~\eqref{dual-advection-diffusion}.}
  \end{figure}

  \cref{fig:primal-advection-diffusion,fig:dual-advection-diffusion} report the evolution of the residual and error on~\eqref{eq:Ax=b} and~\eqref{eq:adjoint-system} for~\eqref{primal-advection-diffusion} and~\eqref{dual-advection-diffusion}, respectively.
  In this numerical illustration, residuals and errors are computed explicitly at each iteration as \(b - Ax\), \(c - A^T t\), \(x - x_\star\), and \(t - t_\star\) in order to discount errors in the approximation formulae for those expressions.
  In this example, \BiLQR terminates in about four times fewer iterations than \TriLQR and six times fewer iterations than \MINRESQLP.
  Only the \USYMLQ error and the \USYMQR residual are monotonic.
  Although the \MINRESQLP residual on~\eqref{eq:adjoint_systems} is monotonic, individual residuals on~\eqref{eq:Ax=b} and~\eqref{eq:adjoint-system} are not.

  We use bicubic spline interpolation and $3 \times 3$ points Gauss quadrature to computate estimates of $J(u)$ with and without correction term.
  With the $u_{\star}$ given above, $J_{\star} := J(u_\star) = (\pi (e + 1))^2 / (\pi^2 + 1)^2$.
  The right plot of \cref{fig:J} illustrates the error in the evaluation of \(J(u)\) as a function of \(h\) using the naive \(J(u) \approx J(u_h)\) and improved~\eqref{eq:superconvergent-approx} approximations.

  \section{Discussion}

  \BiLQ completes the family of Krylov methods based on the Lanczos biorthogonalization process, and is a natural companion to \BiCG and \QMR.
  It is a quasi-minimum error method, and in general, neither the error not the residual norm are monotonic.

  Contrary to the \cite{arnoldi-1951} and the \cite{golub-kahan-1965} processes, the Lanczos biorthogonalization and orthogonal trigonalization processes require two initial vectors.
  This distinguishing feature makes them readily suited to the simultaneous solution of primal and adjoint systems.
  A prime application is the superconvergent estimation of integral functionals in the context of discretized ODEs and PDEs.
  In our experiments, we observed that \BiLQR outperforms both \TriLQR and \MINRESQLP applied to an augmented system in terms of error and residual norms.

  Our Julia implementation of \BiLQ, \QMR, \BiLQR, \TriLQR and \MINRESQLP are available from \https{github.com/JuliaSmoothOptimizers/Krylov.jl} and can be applied in any floating-point arithmetic supported by the language.
  In our experiments with adjoint systems, we run both the primal and ajoint solvers until both residuals are small.
  A slightly more sophisticated implementation would interrupt the first solver that converges and only apply the other until it too converges.
  That is the strategy applied by \cite{buttari-orban-ruiz-titley_peloquin-2019}.

  \MINRES applied to~\eqref{eq:adjoint_systems} does not produce monotonic residuals in the individual primal and adjoint systems.
  In our experiments, we explicitly computed those residuals but \cite{herzog-soodhalter-2017} devised a modification of \MINRES that allows to monitor block residuals that could be of use in the context of estimating integral functionals.

  Although the \BiLQ error is not monotonic in the Euclidean norm, it is in the \(U_p U_p^T\)-norm, which is not iteration dependent, but is unknown until the end of the biorthogonalization process.
  The same property holds for the \QMR residual.
  Exploiting such properties to obtain useful bounds on the \BiLQ and \BiCG error in Euclidean norm that could help devise useful stopping criteria is the subject of ongoing research.

  \small
  \bibliographystyle{abbrvnat}
  \bibliography{abbrv,bilq}
  \normalsize

\end{document}